\documentclass[12pt]{article}

\usepackage{amsmath, amsfonts, amssymb}
\usepackage{latexsym}
\usepackage{graphicx}
\usepackage{color}
\usepackage{url}

\newtheorem{thm}{Theorem}[section]

 \newtheorem{lem}[thm]{Lemma}
 \newtheorem{prop}[thm]{Proposition}
 \newtheorem{exa}[thm]{Example}

\newcommand{\R}{\mathbb{R}}

\newcommand{\bh}{{\mathbf h}}

\newcommand{\bg}{{\mathbf g}}
\newcommand{\bx}{{\mathbf x}}
\newcommand{\bk}{{\mathbf k}}

\newcommand{\bob}{{\mathbf b}}

\newcommand{\bZ}{{\mathbf Z}}

\newcommand{\bR}{{\mathbf R}}

\newcommand{\cA}{{\mathcal A}}

\newcommand{\cB}{{\mathcal B}}

\newcommand{\cL}{{\mathcal L}}

\newcommand{\cS}{{\mathcal S}}
\newcommand{\cO}{{\mathcal O}}

\newcommand{\al}{\alpha}
\newcommand{\be}{\beta}

\newcommand{\ph}{\varphi}

\newcommand{\si}{\sigma}

\newcommand{\n}{\|}

\newcommand{\ti}{\widetilde}

\def\s{{\mathcal S}}
\def\Aut{\operatorname{Aut}}
\def\GL{\operatorname{GL}}

\renewcommand{\thefootnote}{\fnsymbol{footnote}}

\begin{document}

\begin{center}
%{\Large \bf Elementary proofs of some theorems on lattices\\[0.5ex] generated by finite Abelian groups}
%{\Large \bf The Cauchy-Binet formula for determinants and\\[0.5ex] lattices generated by finite Abelian groups}
{\Large {\bf On lattices generated by finite Abelian groups}}

\vspace{7mm}
{\large Albrecht B\"ottcher, Lenny Fukshansky,\\[0.5ex]

Stephan Ramon Garcia, Hiren Maharaj}
\end{center}

\medskip
\begin{quote}
\footnotesize{This paper is devoted to the study of lattices generated by finite
Abelian groups. Special species of such lattices arise in the exploration of elliptic curves over finite fields. In case
the generating group is cyclic, they are also known as the Barnes lattices. It is shown that for every
finite Abelian group with the exception of the cyclic group of order four
these lattices have a basis of minimal vectors. Another result provides an improvement of
a recent upper bound by Min Sha for the covering radius in the case of the Barnes lattices. Also discussed are properties of the automorphism groups of these lattices.
}
\let\thefootnote\relax\footnote{\hspace*{-7.5mm} MSC 2010: Primary 11H31; Secondary 11G20, 11H55, 15A15, 15B05, 52C17}
\let\thefootnote\relax\footnote{\hspace*{-7.5mm} Keywords: well-rounded lattice, finite Abelian group, minimal vector, covering radius,
automorphism group,
Toeplitz determinant}
\let\thefootnote\relax\footnote{\hspace*{-7.5mm}  Fukshansky acknowledges support by Simons Foundation grant \#279155,
Garcia acknowledges support by NSF grant DMS-1265973.}
\end{quote}

\section{Introduction}\label{SI}

The lattice generated by a finite Abelian (additive) group $G=\{0,g_1, \ldots, g_n\}$
of order $|G|=n+1$
is defined as
\[\cL(G):=\{X=(x_1,\ldots, x_n, -x_1-\cdots-x_n)\in \bZ^{n+1}: x_1g_1+\cdots+x_n g_n=0\}.\]
We think of this lattice as a sublattice of full rank $n$ of the root lattice
\[\cA_{n}:=\{(x_1, \ldots, x_n, -x_1-\cdots-x_n)\in \bZ^{n+1}\}.\]
We denote by $d(G)$ the minimum distance in $\cL(G)$, that is, with $\n \cdot \n$
denoting the Euclidean norm,
\[d(G):=\min \{\n X\n: X \in \cL(G)\setminus\{0\}\},\]
and we let $\cS(G)$ stand for the set of nonzero lattice vectors of minimal length,
that is, for the set of all $X \in \cL(G)$ with $\n X\n =d(G)$.
The lattice $\cL(G)$ is said
\begin{itemize}
\item to be well-rounded if $\cS(G)$ contains  $n$ linearly independent vectors,
\item to be generated by minimal vectors if ${\rm span}_\bZ \,\cS(G)=\cL(G)$, that is, if each
vector in $\cL(G)$ is a linear combination with integer coefficients of vectors in $\cS(G)$,
\item to have a basis of minimal vectors if $\cS(G)$ contains $n$ vectors such that each lattice vector is a linear combination
with integer coefficients of these $n$ vectors.
\end{itemize}

Clearly, each of these properties implies its predecessor. Lattices in Euclidean spaces satisfying any of the above properties are of importance in extremal lattice theory, discrete geometry, and combinatorics. Such lattices usually have a high degree of symmetry, which allows for some classical discrete optimization problems to be reduced to them (see~\cite{martinet} for detailed information). It is especially interesting when lattices with these properties come from algebraic constructions, hence inheriting additional algebraic structure. For instance, there are well-known lattice constructions from ideals in number fields~\cite{bayer1}, \cite{bayer2}, ideals in polynomial rings~\cite{lub_mic}, and curves over finite fields~\cite{tv} (pp. 578--583). In addition to their intrinsic theoretical value, such lattices also have many applications, for instance in coding theory and cryptography, as described in~\cite{tv} and~\cite{lub_mic}, respectively.
\medskip

Our present construction of lattices from Abelian groups generalizes the special case of a family of lattices coming from elliptic curves over finite fields as in~\cite{tv}, which has recently been investigated in~\cite{FM} and~\cite{Sha}. It is our goal to show that these lattices have some remarkable geometric properties, including those listed above. Here is our first observation.

\begin{thm} \label{Theo 1.2}
Except for the lattice $\cL(\bZ_4)$, which is not well-rounded, the
lattice $\cL(G)$ is well-rounded for every finite Abelian group $G$. The minimum distance is $\sqrt{8}$ for $G=\bZ_2$, is $\sqrt{6}$
for $G=\bZ_3$, and equals $\sqrt{4}=2$ for all other finite Abelian groups $G$.
\end{thm}

Our first main result is as follows.

\begin{thm} \label{Theo 1.3}
For every finite Abelian group $G \neq \bZ_4$,  the
lattice $\cL(G)$ has a basis of minimal vectors.
\end{thm}

Theorem \ref{Theo 1.2} implies that $\cL(\bZ_4)$ does not possess a basis of minimal vectors.
Of course, Theorem~\ref{Theo 1.3} is stronger than Theorem~\ref{Theo 1.2}. We nevertheless
give an independent proof for Theorem~\ref{Theo 1.2}, because well-roundedness may be proved
by arguments that are much simpler than those we have to invoke to establish Theorem~\ref{Theo 1.3}.

\medskip
One of the subtleties of lattices, discovered
by Conway and Sloane~\cite{CoSlwithout}, is that a lattice generated by minimal vectors need not
have a basis of minimal vectors. More recently, it has been shown~\cite{SM} that this phenomenon takes place for some lattices in dimensions $\geq 10$, but not in lower dimensions. Theorem~\ref{Theo 1.3} implies that this does not happen for the class of lattices explored in this paper.

\medskip
The lattices we study here include the lattices which come from elliptic curves over finite fields, namely $\cL(G)$ where $G$ is the group of rational points on an elliptic curve over a finite field.
These groups were completely described by R\"uck~\cite{Ruck}, and they are always of the form  $G=\bZ_{m_1}\times \bZ_{m_2}$, the direct product of two cyclic groups (with further restrictions
on the possible values of $(m_1,m_2)$).
For lattices coming from elliptic curves over
finite fields, paper~\cite{FM} contains Theorems~\ref{Theo 1.2} and the weaker version of Theorem~\ref{Theo 1.3}
which states that for $|G| \ge 5$ the lattice $\cL(G)$ is generated by minimal vectors, while Sha~\cite{Sha}
proved Theorem~\ref{Theo 1.3} for those lattices. The contribution of the present paper is that we extend these results to
general finite Abelian groups $G$.
\medskip

Well-rounded lattices play a crucial role in the theory of sphere packing (see~\cite{CoSl}, \cite{martinet}), where maximal non-overlapping balls of equal radius (equal to half of the minimal distance of the lattice) are centered at the lattice points with the goal of covering the largest possible proportion of the ambient space. This proportion, called the packing density of the lattice, is equal to the volume of one such ball divided by the volume of a fundamental domain of the lattice (equal to the determinant of the lattice). The lattice packing problem then aims to maximize the packing density on the space of lattices in a given dimension, which emphasizes the importance of knowing minimal distance and determinant of the lattice.

\medskip
Our second topic of investigation is related to another classical optimization problem on lattices, the sphere covering problem (again, see~\cite{CoSl}, \cite{martinet}). The goal is to cover the ambient space completely by balls of equal radius (called the covering radius of the lattice) centered at the lattice points, minimizing the proportion of overlap of these balls. A variety of classical general bounds for covering radii of lattices (also referred to as inhomogeneous minima) can be found in~\cite[Chap. 2, Sec. 13]{gruber}. Here we present estimates for the covering radius $\mu(G)$ of $\cL(G)$.
By definition, $\mu(G)$ is smallest number $\mu$
such that
\[
{\rm span}_{\bR}\cA_n:=\{(\xi_1,\ldots, \xi_{n}, -\xi_1-\cdots-\xi_{n}) \in \bR^{n+1}: \xi_1, \ldots, \xi_n \in \bR\}
\]
is covered by $n$-dimensional closed Euclidean balls of radius $\mu$ centered at the points of $\cL(G)$.
The covering radii for the small groups are
\begin{eqnarray*}
& & \mu(\bZ_2)=\sqrt{2} \approx 1.4142,\\
& & \mu(\bZ_3)=\sqrt{2} \approx 1.4142,\\
& & \mu(\bZ_4)=\frac{3}{2}=1.5000, \quad \mu(\bZ_2 \times \bZ_2)=\sqrt{3}\approx 1.7321,\\
& & \mu(\bZ_5)=\sqrt{2} \approx 1.4142,\\
& & \mu(\bZ_6)=\sqrt{\frac{17}{8}} \approx 1.4577.
\end{eqnarray*}
The distance of the point $(1/2, 1/2, \ldots, 1/2,-n/2)\in {\rm span}_{\bR}\cA_n$
to the lattice $\cA_n$ and thus all the more to the sublattice $\cL(G)$ is at least
\[\frac{1}{2^2}+\frac{1}{2^2}+\cdots+\frac{1}{2^2}=\frac{n}{4},\]
which implies that $\mu(G) \ge (1/2)\,\sqrt{n}$ for every $G$. Actually a little more can be said.
Namely, we obviously have $\mu(G) \ge \mu(\cA_n)$, where $\mu(\cA_n)$ is the covering radius of $\cA_n$,
which is known to be
\[\mu(\cA_n)=\left\{\begin{array}{ll} \frac{1}{2}\,\sqrt{n+1} & \mbox{if} \;\: n \;\:\mbox{is odd},\\[0.5ex]
 \frac{1}{2}\,\sqrt{n+1-1/(n+1)} & \mbox{if} \;\: n \;\:\mbox{is even}; \end{array}\right.\]
see \cite[Chap. 4, Sec. 6.1]{CoSl}.
In \cite{Sha}, it is shown that $\mu(G) \le \mu(\cA_n)+\sqrt{2}$.

\medskip
If $G=\bZ_{n+1}$ is the cyclic group of the numbers $0,1,\ldots,n$ with addition modulo $n+1$, then $\cL(G)$
is the sublattice of $\cA_n$ formed by the points satisfying $x_1+2x_2+\cdots+nx_n=0$ modulo $n+1$. These lattices
probably first appeared in~\cite{Ba} and are therefore frequently referred to as the Barnes lattices.
Here is another main result of this paper. It provides us with an improvement
of the upper bound  $\mu(\cA_n)+\sqrt{2}$ for cyclic groups, that is, for the Barnes lattices.

\begin{thm} \label{Theo 1.4}
For every $n \ge 2$,
\[\mu(\bZ_{n+1}) < \frac{1}{2}\,\sqrt{n+4\log(n-1)+7-4\log 2+10/n}.\]
\end{thm}

\medskip
The data (chopped after the fourth digit after the decimal point)
for several values of $n$ are shown in the table.

\begin{figure}[h]
\[\begin{array}{|r|r|r|r|}
\hline
n & \mu(\cA_n) & \mbox{Theorem}\;\: \ref{Theo 1.4} & \mu(\cA_n)+\sqrt{2}\\
\hline
 3  & 1.0000  &  1.8257  &  2.4142\\
  4  &  1.0954  &  1.9443 &   2.5097\\
   5  &  1.2247  &  2.0477   & 2.6390\\
    6  &  1.3093  &  2.1408   & 2.7235\\
     20  &  2.2887  &  3.0210  &  3.7029\\
      50   & 3.5700  &  4.1831  &  4.9842\\
      100   & 5.0247  &  5.5387  &  6.4389\\
      1\,000      &  15.8193 & 16.0613  & 17.2335\\
      10\,000      &  50.0025 &  50.1026  & 51.4167\\
      100\,000      &  158.1147 & 158.1536 &  159.5289\\
      1\,000\,000    &  500.0002 & 500.0149 & 501.4145\\
      \hline
\end{array}\]
%\centerline{Table 1}
\end{figure}

\medskip
Thirdly, we investigate a certain property of the automorphism groups of our lattices $\cL(G)$, which is intrinsically related to their algebraic construction.
The automorphism group $\Aut(\cL)$ of a full rank sublattice $\cL$ of some lattice $\cA$ is defined as the group of all maps of $\cL$ onto itself
which extend to linear isometries
of ${\rm span}_\bR \cA$. It is easily seen that in our setting, $\cL(G) \subset \cA_n$, a map $\tau \in \Aut(\cL(G))$
is necessarily of the form
\[\tau(X)=\tau\left(x_1, \ldots, x_n, -\sum_{i=1}^n x_i\right)=\left(Ux,-\sum_{i=1}^n (Ux)_i\right)\]
with some matrix $U \in \GL_n(\bZ)$.  We therefore identify $\Aut(\cL(G))$ with a subgroup of $\GL_n(\bZ)$.
It is a well known fact that any finite subgroup of $\GL_n(\bZ)$ is the automorphism group of some lattice. In all dimensions except for $n=2,4,6,7,8,9,10$ (dimensions with exceptionally symmetric lattices) the largest such group is $(\bZ/2\bZ)^n \rtimes S_n$, the automorphism group of the integer lattice~$\bZ^n$; here $S_n$ is the symmetric group on $n$ letters viewed as the subgroup of $\GL_n(\bZ)$ consisting of the permutation matrices (see~\cite{CoSl}, \cite{martinet}, and~\cite{schur_book} for more information on automorphism groups of lattices). Lattices with large automorphism groups usually have a large degree of geometric symmetry, which often correlates with having many minimal vectors and well-roundedness. In particular, the relation between certain properties of $\Aut(\cL) \cap S_n$ and the probability of $\cL$ being well-rounded has recently been investigated in~\cite{FS-1}, \cite{FS-2}. Here we prove the following.

\begin{thm} \label{Theo 1.5}
For every finite Abelian group $G$,
\[ \Aut(\cL(G)) \cap S_n \cong \Aut(G),\]
where $\Aut(G)$ is the group of automorphisms of $G$.
\end{thm}

This result, along with a characterization of the automorphism groups of finite Abelian groups,
for which see, e.g., \cite{hillar}, helps to understand the symmetries of our family of lattices~$\cL(G)$.

\medskip
The paper is organized as follows.
Section \ref{SDet} is devoted to the determinant of $\cL(G)$. There we first present a short
derivation based on a general fact from lattice theory and then give a second proof,
which uses only elementary facts for determinants, mainly the Cauchy-Binet formula.
The proofs of Theorems~\ref{Theo 1.2} and~\ref{Theo 1.3} we give here occupy Sections~\ref{SP}
to~\ref{SMin} and use tools from linear algebra only. Again the Cauchy-Binet formula is always the key.
Theorem~\ref{Theo 1.4} is proved in Section~\ref{SCov}. The proof is anew pure linear algebra and makes use of explicit formulas for certain Toeplitz determinants. Finally, in Section~\ref{aut} we prove Theorem~\ref{Theo 1.5} and comment on a certain geometric interpretation of this result.

\section{The determinant} \label{SDet}

A set of $n$ vectors $X_1, \ldots, X_n \in \cL(G)$ is called a basis
if each vector in  $\cL(G)$ is a linear combination with integer coefficients of these vectors.
In that case the parallelotope spanned by $X_1, \ldots, X_n$ is referred to as a fundamental
parallelotope. All fundamental parallelotopes have the same volume. This volume
is denoted by $\det \cL(G)$ and
referred to as the determinant of the lattice $\cL(G)$.
Even more can be said:  the parallelotope spanned by $n$ vectors $X_1, \ldots, X_n \in \cL(G)$
has the volume $\det \cL(G)$ if {\em and only if} these vectors form a basis of $\cL(G)$.
If $X_1, \ldots, X_n \in \cL(G)$ form a basis, then the $(n+1) \times n$ matrix $B$ whose
$j$th column is constituted by the $n+1$
coordinates of $X_j$ is called a basis matrix. If
$B$ is an arbitrary basis matrix of $\cL(G)$, then $\det \cL(G)=\sqrt{\det B^\top B}$,
where the determinant on the right is the usual determinant of an $n \times n$ matrix.
All these results are standard in lattice theory and can be found in~\cite{CoSl}, \cite{gruber}, or~\cite{martinet}, for example.

\medskip
It turns out that $\det \cL(G)=|G|^{3/2}=(n+1)^{3/2}$.
The following proof of this formula is from~\cite{Sha}(Proposition 5.1), where it is given for the case when $G$ is a subgroup
of the group of rational points on an elliptic curve over a finite field; it also holds
verbatim for general Abelian groups $G$.
Let $\cL$ be a sublattice of full rank of some lattice $\cA \subset \bR^N$. Think of $\cA$ as an (Abelian)
additive group and consider $\cL$ as a subgroup of $\cA$.
A basic result of lattice theory says that if the quotient group $\cA/\cL$ has finite order $|\cA/\cL|$,
then $\det \cL/\det\cA=|\cA/\cL|$.
Now take $\cA=\cA_n$ and $\cL=\cL(G)$. It is known that $\det \cA_n=\sqrt{n+1}$. The group
homomorphism
\[\ph: \cA_n \to G, \quad (x_1, \ldots, x_n, -x_1-\cdots-x_n) \mapsto x_1g_1+\cdots x_n g_n\]
is surjective and its kernel is just $\cL(G)$. Consequently, $\cA_n/\cL(G)$ is isomorphic to $G$,
which implies that $|\cA_n/\cL(G)|=|G|=n+1$. It follows that
\[\frac{\det \cL(G)}{\sqrt{n+1}}=\frac{\det\cL(G)}{\det\cA_n}=|\cA_n/\cL(G)|=n+1,\]
as asserted.

\medskip
Here is a purely linear algebra proof of the same determinant formula. We first exemplify the idea by considering $G=\bZ_2 \times \bZ_4$.
The lattice $\cL(G)$ consists of the points
\[(x_1,x_2,y_{02},y_{03},y_{11}, y_{12}, y_{13},-x_1-x_2-y_{02}-y_{03}-y_{11}- y_{12}- y_{13})
\in \bZ^8\]
satisfying
\[x_1(1,0)+x_2(0,1)+y_{02}(0,2)+y_{03}(0,3)+y_{11}(1,1)+ y_{12}(1,2)+ y_{13}(1,3)=(0_2,0_4),\]
where $0_2$ and $0_4$ are the zeros in $\bZ_2$ and $\bZ_4$. We may choose the five numbers $y_{jk}$ arbitrarily,
after which $x_1$ and $x_2$ are determined uniquely modulo $2$ and $4$, respectively. Taking $y_{jk}=1$
and $y_{\al,\be}=0$ for $(\al,\be)\neq (j,k)$, we get $x_1+j=0_2$ and $x_2+k=0_4$, that is, $x_1=-j$ modulo $2$
and $x_2=-k$ modulo $4$. Thus, a basis in $\cL(G)$ is formed by the five rows
\[(-j,-k,0,\ldots,0,1,0, \ldots,0,j+k-1),\]
the number $1$ being at the $(j,k)$th position in lexicographic order, and by the two rows
\[(2,0,0,0,0,0,0,-2), \quad (0,4,0,0,0,0.-4),\]
which allow us move $x_1$ and $x_2$ within $2\bZ$ and $4\bZ$.
It follows that the matrix $B^\top$ formed by these seven rows,
\[B^\top=\left(\begin{array}{rrrrrrrr}
2 & 0 &  &  &  &  &  &-2\\
0 & 4 &  &  &  &  &  &-4\\
0 & -2 & 1&  &  &  &  &1\\
0 & -3 &  & 1&  &  &  &2\\
-1 & -1 &  &  & 1&  &  &1\\
-1 & -2 &  &  &  &  1& &2\\
-1 & -3 &  &  &  &  & 1&3
\end{array}\right),\]
is the transpose of a basis matrix $B$ of the lattice $\cL(G)$. Thus,
$\det \cL(G)=\sqrt{\det B^\top B}$. The Cauchy-Binet formula
gives
\[\det B^\top B = (\det B_1)^2+(\det B_2)^2 + \sum_{j,k}(\det B_{jk})^2+(\det B_7)^2,\]
where $B_1, B_2,B_7$ result from $B^\top$ by deleting the columns $1,2,7$ and $B_{jk}$ is the
matrix obtained by deleting the column with $1$ in the position $(j,k)$. Expanding the determinants of $B_1,B_2,B_7$
along the five columns with a single $1$, we see that the squares of these determinants are
\[\left|\begin{array}{rr}0 & -2\\4 &-4\end{array}\right|^2=8^2,
\quad \left|\begin{array}{rr}2 & -2\\0 &-4\end{array}\right|^2=8^2,
\quad \left|\begin{array}{rr}2 & -0\\0 &4\end{array}\right|^2=8^2,\]
and expanding the five determinants $\det B_{jk}$ along their four columns with a single $1$ we get
\[(\det B_{jk})^2=\left|\begin{array}{rrr}2 & 0 & -2\\
0 & 4 & -4\\
-j & -k & j+k-1\end{array}\right|^2
= \left|\begin{array}{rrr}2 & 0 & 0\\
0 & 4 & 0\\
-j & -k &-1\end{array}\right|^2=8^2.\]
Thus, $\det B^\top B=8\cdot 8^2=8^3=|G|^3$, as desired.

\medskip
It is clear how to proceed in the general case $G=\bZ_{m_1} \times \cdots \times \bZ_{m_k}$.
Put $m=m_1\cdots m_k$.
Then $B^\top$ has $m-1$ rows and $m$ columns and we may employ the
Cauchy-Binet formula to express $\det B^\top B$ as the sum of $m$ squares of determinants
as above. The first $k$ and the last squared determinants are readily seen to
be $(m_1\cdots m_k)^2=m^2$, and
and the $m-k-1$ squared determinants corresponding to indices $(j_1, \ldots, j_k)$ are,
with $\si:=j_1+\cdots+j_k$,
\[\left|\!\begin{array}{rrrrr}
m_1 & & & &    -m_1\\
 & m_2 &  &       & -m_2\\
 & & \ddots & &  \vdots\\
 & & & m_k & -m_k\\
 -j_1 & -j_2 & \ldots & -j_k & \si-1
\end{array}\!\right|^2
=\left|\!\begin{array}{rrrrr}
m_1 & & & &    0\\
 & m_2 &  &       & 0\\
 & & \ddots & &  \vdots\\
 & & & m_k & 0\\
 -j_1 & -j_2 & \ldots & -j_k & -1
\end{array}\!\right|^2=(m_1\cdots m_k)^2=m^2.
\]
Consequently, $\det B^\top B=m\cdot m^2=m^3=|G|^3$.

\section{The small groups}\label{SP}

We now turn to the proof of Theorems \ref{Theo 1.2} and \ref{Theo 1.3}.

\medskip
We arrange the nonzero elements of $G$ in a column $\bg=(g_1, \ldots, g_n)^\top$ of height $n$.
Obviously, there are $n!$ possibilities to do this. Then
each point \[X=(x_1,\ldots, x_n, -(x_1+\cdots+x_n))\in\cL(G)\] may be represented
by a column $\bx=(x_1, \ldots, x_n)^\top$ of the same height $n$.
Given $n$ points $X_1, \ldots, X_n$ in $\cL(G)$, we denote by $M$ the $n \times n$
matrix composed of the columns $\bx_1, \ldots, \bx_n$, and we collect the data
in an array $\bg||M$. The arrays that may be obtained in this way
will be called admissible for $G$. We let $\ti{M}$ stand for the $(n+1) \times n$ matrix
which results from adding the row consisting of the negatives of the column sums of
the matrix $M$. Thus, $n$ vectors $X_1, \ldots, X_n$ form a basis in $\cL(G)$ if
and only if $\ti{M}$ is a basis matrix, which, because $|G|^3=(\det\cL(G))^2$, is equivalent to the
equality $\det \ti{M}^\top \ti{M} = |G|^3$.

\medskip
Clearly, $n$ vectors $X_1, \ldots, X_n$
are linearly independent if and only if so are the $n$ columns $\bx_1, \ldots, \bx_n$.
Therefore, in order to prove that $\cL(G)$ is well-rounded, we have to find an admissible array $\bg||M$ in which the matrix $M$
is nonsingular and comes from points of minimum distance. To prove the stronger property that
$\cL(G)$ has a basis of minimal vectors, we have to find an admissible array $\bg||M$
associated with points of minimum distance such that $\ti{M}$ is a basis matrix.

\medskip
First of all we remark that the minimum distance is always at least $\sqrt{1^2\!+\!1^2\!+\!1^2\!+\!1^2}=2$
and that this distance is attained exactly at the points $X$ containing two times $1$, two times $-1$,
and otherwise only zeros. If the lattice does not contain such points, the minimum distance must
be at least $\sqrt{1^2+1^2+2^2}=\sqrt{6}$.

\begin{exa} \label{Exa 2.1}{\rm
Let $G$ be $\bZ_3=\{0,1,2\}$ and let $\bg=(1,2)^\top$. (The other possibility would be to put
$\bg=(2,1)^\top$.) Then $\cL(\bZ_3)$
consists of the integer points $(x,y,-x-y)$ satisfying $x+2y=0$ modulo $3$. By inspection it is easily seen
that $d(\bZ_3)$ is $\sqrt{6}$ and that exactly six points of $\cL(\bZ_3)$ have minimal distance.
Two of them are the points $X_1=(-2,1,1)$ and $X_2=(1,-2,1)$. The
array $\bg||M$ corresponding the these two points is
\[\begin{array}{r||rr}
1 & -2 & 1\\
2 & 1 & -2
\end{array}.\]
The matrix $M=\big(\!\!\begin{array}{rr} -2 & \!\!1\\[-1ex] 1 & \!\!-2 \end{array}\!\!\big)$
is nonsingular, and hence $\cL(\bZ_3)$ is well-rounded with the minimum distance $\sqrt{6}$.
We have
\[\ti{M}=\left(\begin{array}{rr}
-2 & 1\\
1 & -2\\
\hline
1 & 1\end{array}\right).\]
Since $\det\ti{M}^\top \ti{M}=3^3$, we see that $\ti{M}$ is a basis matrix and thus that
$\cL(\bZ_3)$ has a basis of minimal vectors.

}\end{exa}

\begin{exa} \label{Exa 2.2}{\rm
Things are trivial for $G=\bZ_2$, in which case $n=1$. We have $\cL(\bZ_2)=\{(2x,-2x): x \in \bZ\}$,
the minimum distance is $d(\bZ_2)=\sqrt{2^2+2^2}=\sqrt{8}$, and it is attained for $X=(-2,2)$
(and also for $X=(2,-2)$).
}\end{exa}

\begin{exa} \label{Exa 2.3}{\rm
Let $G=\bZ_4$ and $\bg=(1,2,3)^\top$. An integer point $(x,y,z,-x-y-z)$ is in $\cL(\bZ_4)$
if and only if $x+2y+3z=0$ modulo $4$. The points of minimum distance are
\[X_1=(1,1,-1,-1), \;\: X_2=(-1,1,1,-1), \;\: X_3=(-1,-1,1,1), \;\: X_4=(1,-1,-1,1),\]
but any three of them are linearly dependent.
Thus, $\cL(\bZ_4)$ is not well-rounded. Clearly, $d(\bZ_4)=2$.
}\end{exa}

\begin{exa} \label{Exa 2.4}{\rm
For $G=\bZ_2 \times \bZ_2$, the array
\[\bg||M = \begin{array}{r||rrr}
(0,1) & 1 & -1 & 1\\
(1,0) & 1 & 1 & -1\\
(1,1) & -1 & 1 & 1\end{array}\]
is admissible, and since $\det M=4 \neq 0$, it follows that $\cL(\bZ_2 \times \bZ_2)$ is well-rounded
with $d(\bZ_2 \times \bZ_2)=2$. The matrix
\[\ti{M}=\left(\begin{array}{rrr}
1 & -1 & 1\\
1 & 1 & -1\\
-1 & 1 & 1\\
\hline
-1 & -1 & -1\end{array}\right)\]
satisfies $\det \ti{M}^\top \ti{M}=4^3$, and hence $\cL(\bZ_2 \times \bZ_2)$ has
a basis of minimal vectors.
}\end{exa}

\section{The cyclic groups}\label{SC}

Let $G=\bZ_{m}$ with $m \ge 5$ and put $\bg_m=(1,2,\ldots,m-1)^\top$.
We denote by $T_m$ the
$(m-1) \times (m-1)$ tridiagonal Toeplitz matrix with $-2$ on the main diagonal
and $1$ on the two neighboring diagonals. For example,
\[T_7=\left(\begin{array}{rrrrrr}
-2 & 1 & 0 & 0 & 0 & 0\\
1 & -2 & 1 & 0 & 0 & 0\\
0 & 1 & -2 & 1 & 0 & 0\\
0 & 0 & 1 & -2 & 1 & 0\\
0 & 0 & 0 & 1 & -2 & 1\\
0 & 0 & 0 & 0 & 1 & -2\end{array}\right),\quad
\ti{T}_7=\left(\begin{array}{rrrrrr}
-2 & 1 & 0 & 0 & 0 & 0\\
1 & -2 & 1 & 0 & 0 & 0\\
0 & 1 & -2 & 1 & 0 & 0\\
0 & 0 & 1 & -2 & 1 & 0\\
0 & 0 & 0 & 1 & -2 & 1\\
0 & 0 & 0 & 0 & 1 & -2\\
\hline
1 & 0 & 0 & 0 & 0 & 1\end{array}\right).\]
Let $U_m$ be the $(m-1) \times (m-1)$ matrix which results from the $(m-1) \times (m-1)$ bidiagonal Toeplitz
matrix with $1$ on the main diagonal and on the subdiagonal after replacing the last column with
$(0,\ldots,0,-1,-1,-1,0)^\top$. For instance,
\[U_5=\left(\begin{array}{rrrr}
1 & 0 & 0 & -1\\
1 & 1 & 0 & -1\\
0 & 1 & 1 & -1\\
0 & 0 & 1 & 0
\end{array}\right),\quad
U_7=\left(\begin{array}{rr|rrrr}
1 & 0 & 0 &  0 & 0 & 0\\
1 & 1 & 0 & 0 & 0 & 0\\
\hline
0 & 1 & 1 & 0 & 0 &-1\\
0 & 0 & 1 & 1 & 0 & -1\\
0 & 0 & 0 & 1 & 1 & -1\\
0 & 0 & 0 & 0 & 1 & 0
\end{array}\right).\]

\begin{lem}\label{Lem Cyc1} Let
$M_m=T_mU_m$ and $\ti{M}_m=\ti{T}_mU_m$.
Then the matrix $M_m$ results from the $(m-1) \times (m-1)$ tetradiagonal Toeplitz
matrix with first column $(-1,-1,1, 0,\ldots,0)^\top$ and first row $(-1,1,0,\ldots,0)$ by replacing the
last column with $(1,0,1,-1)^\top$ for $m=5$ and with the column $(0, \ldots,0,-1,1,0,1,-1)^\top$ for $m \ge 6$.
\end{lem}

{\em Proof.} Direct computation. $\;\:\square$

\medskip
It can be checked straightforwardly that $\bg_m||M_m$
is an admissible array for $\bZ_m$. For example, the arrays $\bg_5||M_5$ and $\bg_7||M_7$ are
\[\bg_5||M_5=\begin{array}{r||rrrr}
1 & -1 & 1 & 0 & 1\\
2 & -1 & -1 & 1 & 0\\
3 & 1 & -1 & -1 & 1\\
4 & 0 & 1 & -1 & -1 \end{array},
\quad
\bg_7||M_7=\begin{array}{r||rrrrrr}
1 & -1 & 1 & 0 & 0 & 0 & 0 \\
2 & -1 & -1 & 1 & 0 & 0 & -1\\
3 & 1 & -1 & -1 & 1& 0 & 1\\
4 & 0  & 1 & -1 & -1 & 1 & 0\\
5 & 0 & 0 & 1 & -1 & -1 & 1\\
6 & 0 & 0 & 0& 1 & -1 & -1\end{array}.\]
It is well known that $\det T_m=(-1)^{m-1}m$. We have $\det U_5=1$,
which implies that $\det U_m=1$ for all $m \ge 5$. Consequently, by Lemma~\ref{Lem Cyc1},
\[\det M_m =\det T_m \det U_m =(-1)^{m-1}m \neq 0.\]
This proves that $\cL(\bZ_m)$ is well-rounded with $d(\bZ_m)=2$. The fact that $\cL(\bZ_m)$
has a basis of minimal vectors lies a little deeper. It requires the following result.

\begin{lem} \label{Lem Cyc2}
We have $\det\ti{T}_m^\top \ti{T}_m=m^3$.
\end{lem}

{\em Proof.}
Applying the Cauchy-Binet formula, we may write
\[\det \ti{T}_m^\top \ti{T}_m =(\det C_1)^2+(\det C_2)^2+\cdots+(\det C_m)^2,\]
where $C_j$ results from $\ti{T}_m$ by deleting the $j$th row.
Clearly, $(\det C_m)^2=(\det T_m)^2=m^2$. For $j \le m-1$, we
expand $\det C_j$ along the last row and obtain two block-triangular determinants:
\[\det C_j=(-1)^m\det T_{m-j}+\det T_j=(-1)^m(-1)^{m-j-1}(m-j)+(-1)^{j-1}j=(-1)^{j-1}m.\]
It follows that $(\det C_j)^2=m^2$.
Consequently, $\det  \ti{T}_m^\top \ti{T}_m = m\cdot m^2=m^3$.
$\;\:\square$

\medskip
Combining Lemma \ref{Lem Cyc2} with the factorization $\ti{M_m}=\ti{T_m}U_m$ delivered by Lemma~\ref{Lem Cyc1},
we get
$\det\ti{M}_m^\top\ti{M}_m=m^3$, which shows that $\bZ_m$ is generated by vectors of minimum distance.

\section{Direct products: well-roundedness} \label{SD}

In this section we complete the proof of Theorem \ref{Theo 1.2}. Much of the following, especially the choice
of the matrices in the arrays, resembles the constructions in~\cite{Sha}. However, our reasoning is consistently based on
the computation of determinants and thus completely differs from the arguments used in~\cite{Sha}.

\begin{lem} \label{Lem 4.1}
If $G$ and $H$ are finite Abelian groups such that $\cL(G)$ and $\cL(H)$ are well-rounded
with $d(G)=d(H)=2$, then $\cL(G \times H)$ is well-rounded and $d(G\times H)=2$.
\end{lem}

{\em Proof.} Let $G=\{0,g_1, \ldots, g_n\}$ and $H=\{0,h_1, \ldots, h_m\}$.
We write $\bg=(g_1, \ldots,g_n)^\top$ and $\bh=(h_1,\ldots,h_m)^\top$.
By assumption,
there exist nonsingular integer matrices $M_G=(a_{ij})$ and $M_H=(b_{ij})$ such that $\bg||M_G$ and
$\bh||M_H$ are admissible arrays and such that the columns after deleting all zeros
reduce to columns with $3$ or $4$ entries containing only $\pm1$ and having their column sum in $\{-1,0,1\}$. The array
\[\begin{array}{l||rrr|rrr|rrrr}
(g_1,0) & a_{11} & \ldots & a_{1n} & & & & -1 & -1 &\\
\vdots & \vdots & & \vdots & & & & & &\\
(g_n,0) & a_{n1} & \ldots & a_{nn}& & & & & & \vdots\\
\hline
(0,h_1) & & & & b_{11} & \ldots & b_{1m} & -1 & & \\
(0,h_2) & & & & b_{21} & \ldots & b_{2m} &  & -1 &\\
\vdots & & & & \vdots & & \vdots & & & \\
(0,h_m) & & & & b_{m1} & \ldots & b_{mm} & & & \vdots \\
\hline
(g_1,h_1) & & & & & & & 1 & & \\
(g_1,h_2) & & & & & & & & 1 & \\
\vdots & & & & & & & & & \vdots
\end{array}
\]
consists of  $n+m+nm=(n+1)(m+1)-1$ columns. The last $nm$ columns may be labelled by $(g_i,h_j)$,
and the column with this label has $1$ at position $(g_i,h_j)$
and $-1$ at the positions $(g_i,0)$ and $(0,h_j)$.
This array is clearly admissible for $G \times H$, and since its matrix is upper block-triangular
with determinant $\det M_G \det M_H \neq 0$, we conclude that $\cL(G \times H)$ is
well-rounded. This array also reveals that $d(G \times H)=2$. $\;\:\square$

\medskip
Lemma \ref{Lem 4.1} in conjunction with the result of the previous section
proves Theorem~\ref{Theo 1.2} for all groups which do not contain the factors
$\bZ_2, \bZ_3, \bZ_4$.

\begin{lem} \label{Lem 4.2}
If $m \in \{2,3,4\}$ and $G$ is a finite Abelian group such that $\cL(G)$ is well-rounded
with $d(G)=2$, then $\cL(\bZ_m \times G)$ is well-rounded and $d(\bZ_m\times G)=2$.
\end{lem}

{\em Proof.} Let $G=\{0,g_1, \ldots, g_n\}$ and $\bg=(g_1, \ldots, g_n)^\top$.
By the examples in Section~\ref{SP}, we may assume that $n \ge 3$.
Take an admissible array $\bg||M$ with a nonsingular integer matrix $M=(a_{ij})$.
The columns of $M$ may be assumed to be as
described in the preceding proof. The array
\[\begin{array}{l||rrr|rr|rrr}
(0,g_1) & a_{11} & \ldots  & a_{1n} & -1 & -1 &  &  & \\
(0,g_2) & a_{21} & & a_{2n}& & & -1 & & \\
(0,g_3) & a_{31} & & a_{3n} & & & & -1 & \\
\vdots & \vdots & & \vdots & & & & & \\
(0,g_n) & a_{n1} & \ldots  & a_{nn}& & & & & \vdots \\
\hline
(1,0)& & & & 1 & -1 & -1 & -1 & \\
(1,g_1) & & & & 1 & 1 & & & \\
\hline
(1,g_2) & & & &  &  & 1 &  & \\
(1,g_3) & & & &  & & & 1 & \\
\vdots & & & & & & & & \vdots \\
\end{array}\]
is admissible for $\bZ_2 \times G$. The matrix is upper-block triangular with
determinant $\det M \cdot 2 \neq 0$, and we have $d(\bZ_2 \times G)=2$. We turn to $\bZ_3 \times G$.
Suppose $g_1+g_2=0$ and $g_1+g_1=g_3$. Then the array
\[\begin{array}{l||rrr|rrr|rrr|rrrr}
(0,g_1)& a_{11}& \ldots & a_{1n} & &   &  -1    & & & & & & &\\
(0,g_2)& a_{21} & & a_{2n}       & &    &    & -1 & & &-1 &  & & \\
(0,g_3)& a_{31} & & a_{3n}       & & -1   &    &     & \vdots & & & & & \\
\vdots & \vdots  & & \vdots      & &    &     &    &        &  & &-1 &  &\\
(0,g_n) & a_{n1} & \ldots & a_{nn}      &     &    &   &  & & -1 & & & & \vdots\\
\hline
(1,0) & & & & 1 & 0     & -1 & -1 & & -1 &  & & &\\
(1,g_1) & & & & 1 & 1 & 1  &    & \vdots & &  & & &\\
(2,g_1) & & & & -1 & 1 & 0 &    & &  & -1 & -1 & &\\
\hline
(1,g_2) & & & &  & & & 1 & & &  & & & \\
\vdots  & & & & & & & &   \vdots & & & & & \\
(1,g_n) & & & & & & & &        & 1 & & & & \vdots \\
\hline
(2,0) & & & & & & & & & & 1 &   & &\\
(2,g_2) & & & & & & & & & & & 1 & &\\
\vdots & & & & & & & & & & &  & & \vdots \\
\end{array}
\]
is admissible. The matrix of this array is upper block-triangular. The determinant of the $n \times n$ block is nonzero,
and the determinant of the $3\times 3$ block equals $-3$. Thus, $\cL(\bZ_3 \times G)$ is well-rounded
with $d(\bZ_3 \times G)=2$. We finally consider $\bZ_4 \times G$. Let $g_1+g_n=0$. Now the array
\[\begin{array}{l||rrr|rrr|rrrrr}
(0,g_1)& a_{11} & \ldots & a_{1n} &    &    & -1 & -1    & -1 & -1 & & \\
(0,g_2)& a_{21} &        & a_{2n} &    &    &    &     & & & -1 & \\
\vdots &        &        & \vdots &    &    &    &     & & & & \\
(0,g_n)& a_{n1} & \ldots & a_{nn} &    &    & -1 &     & & & & \vdots\\
\hline
(1,0)  &        &         &       &  1 & -1 & 1  & -1  & & & -1 & \\
(2,0)  &        &         &       &  1 & 1  & 0  &     &  -1 & & & \\
(3,0)  &        &         &       & -1 & 1  & 1  &     & & -1 &  &  \\
\hline
(1,g_1)&        &         &       &    &    &    & 1   & & & & \\
(2,g_1)& & & & & & & & 1 & & & \\
(3,g_1)& & & & & & & & & 1 & &  \\
(1,g_2)& & & & & & & & & & 1 & \\
\vdots & & & & & & & & & & & \vdots \\
\end{array}\]
is admissible. The determinant of the $3 \times 3$ block is $4$ and thus nonzero. It follows that
$\cL(\bZ_4 \times G)$ is well-rounded with $d(\bZ_4 \times G)=2$. $\:\;\square$

\begin{lem} \label{Lem 4.3}
The lattices $\cL(\bZ_2 \times \bZ_4)$, $\cL(\bZ_3 \times \bZ_3)$, and $\cL(\bZ_4 \times \bZ_4)$ are
well-rounded with minimum distance $2$.
\end{lem}

{\em Proof.} The array
\[\begin{array}{l||rrrrr|rr}
(1,0)& 1& & & & 1 & & -1\\
(1,3)& 1& 1& & & -1& & \\
(0,3)& -1& 1& 1& -1& & -1& \\
(1,2)& & -1& 1& 1& -1& & -1\\
(1,1)& & & -1& 1& 1& & \\
\hline
(0,1)& & & & & & 1& 0\\
(0,2)& & & & & &1 & 1\\
\end{array}\]
is admissible for $\bZ_2 \times \bZ_4$, the determinants of the diagonal blocks
being $8$ and $1$, which proves the assertion for $\bZ_2 \times \bZ_4$. The array
\[\begin{array}{l||rrrrrrrr}
(0,1)& 1 & & & & & & -1&1 \\
(1,0)& 1& 1& & & & & & -1\\
(1,1)& -1& 1& 1& & & & & \\
(2,1)& &-1 & 1& 1& & & & \\
(0,2)& & & -1& 1& 1& & & \\
(2,0)& & & & -1& 1& 1& & \\
(2,2)& & & & &  -1& 1& 1&\\
(1,2)& & & & & & -1 & 1 & 1\\
\end{array}\]
is admissible for $\bZ_3 \times \bZ_3$, and since the determinant of the entire $8 \times 8$
matrix is $-45$, we get the assertion in this case. Finally, the array
\[\begin{array}{l||rrrrrr|rrrrrr|rrr}
(0,1) & 1 &   &  & & -1 & 1 & & & & & & & & & \\
(1,0) & 1 & 1 &  & &    &-1 & & & & & & & & &\\
(1,1) &-1 & 1 &1 & &    &   & & & & & & & & &\\
(2,1) &   &-1 &1 &1&    &   & & & & & & & & & \\
(3,2) &   &   &-1&1&1   &   & & & & & & & & & \\
(1,3) &   &   &  &-1& 1 &1  & & & & & & & & & \\
\hline
(1,2) & & & & & & & 1 &   &  & & -1 & 1 & & &  \\
(3,1) & & & & & & & 1 & 1 &  & &    &-1 & & & \\
(0,3) & & & & & & &-1 & 1 &1 & &    &   & & & \\
(3,0) & & & & & & &   &-1 &1 &1&    &   & & &  \\
(3,3) & & & & & & &   &   &-1&1&1   &   & & &  \\
(2,3) & & & & & & &   &   &  &-1& 1 &1  & & &  \\
\hline
(2,0) & & & & & & & & & &  &   &  &1 & -1 & 1   \\
(0,2) & & & & & & & & & &  &  &  &1 &  1  &-1  \\
(2,2) & & &  & & & & & & & &  & & -1&  1  &1
\end{array}\]
is admissible for $\bZ_4 \times \bZ_4$. The determinants of
the diagonal blocks are $-16$, $-16$, $4$. Consequently, $\cL(\bZ_4 \times \bZ_4)$
is well-rounded with minimum distance $2$. $\;\:\square$

\medskip
Now we can finish the game. Let
\[G=\underbrace{\bZ_2 \times \cdots \times \bZ_2}_{i}\times
\underbrace{\bZ_3 \times \cdots \times \bZ_3}_{j}\times
\underbrace{\bZ_4 \times \cdots \times \bZ_4}_{k}\times H,\]
where $H$ contains only cyclic groups of order at least $5$ or where
$H$ is absent. In the former case repeated application of Lemmas~\ref{Lem 4.1} and~\ref{Lem 4.2}
shows that $\cL(G)$ is well-rounded with $d(G)=2$. We are left with the latter case. Since
$\bZ_2 \times \bZ_3=\bZ_6$, $\bZ_3 \times \bZ_4=\bZ_{12}$, $\bZ_2 \times \bZ_4$ are well-rounded
with minimum distance $2$ (Section~\ref{SC} for the first two and Lemma~\ref{Lem 4.3} for the last group),
Lemmas~\ref{Lem 4.1} and~\ref{Lem 4.2} give the assertion if two of the numbers $i,j,k$ are at least $1$.
It remains to consider the cases where $G$ is one of the goups
\[G_2=\underbrace{\bZ_2 \times \cdots \times \bZ_2}_{i}, \quad
G_3=\underbrace{\bZ_3 \times \cdots \times \bZ_3}_{j}, \quad
G_4=\underbrace{\bZ_4 \times \cdots \times \bZ_4}_{k}.\]
For $i=1$ we are in Example~\ref{Exa 2.2}, and for $i \ge 2$ we obtain from Example~\ref{Exa 2.4} and
Lemma~\ref{Lem 4.2} that $G_2$ is as asserted. The case of $G_3$ is settled by Example~\ref{Exa 2.1}
for $j=1$ and by Lemmas~\ref{Lem 4.2} and~\ref{Lem 4.3} for $j \ge 2$. Example~\ref{Exa 2.3} ($k=1$)
and Lemmas~\ref{Lem 4.2} and~\ref{Lem 4.3} ($k \ge 2$) finally yield the assertion for $G_4$.

\section{Direct products: bases of minimal vectors} \label{SMin}

This section is devoted to the proof of Theorem \ref{Theo 1.3}. We want to emphasize once more
that Theorem~\ref{Theo 1.3} was previously proved by Sha~\cite{Sha} for subgroups $G$ of the direct product of two
cyclic groups. In particular, Lemma~\ref{Lem 4.3a} and results resembling Lemmas~\ref{Lem 4.1a} and~\ref{Lem 4.2a} in the cases of cyclic groups $G$ and $H$
were already established in~\cite{Sha} using arguments different from ours.

\begin{lem} \label{Lem 4.1a}
Let $G$ and $H$ be finite Abelian groups such that $\cL(G)$ and $\cL(H)$ have bases of minimal vectors
and such that $d(G)=d(H)=2$. Also suppose that there are admissible arrays $\bg||M_G$ and $\bh||M_H$
coming from minimal basis vectors such that $\det M_G=\pm |G|$ and $\det M_H=\pm |H|$. Put $K=G \times H$.
Then $\cL(K)$ has a basis of minimal vectors, $d(K)=2$,
and there exists an admissible array $\bk||M_K$ resulting from minimal basis vectors such that $\det M_K=\pm |K|$.
\end{lem}

{\em Proof.} Let $G,H,\bg,\bh$ be as in the proof of Lemma~\ref{Lem 4.1}. Our present assumptions guarantee
that the two matrices $M_G$ and $M_H$ in the proof of Lemma~\ref{Lem 4.1} may be taken so that
$\ti{M}_G$ and $\ti{M}_H$ are basis matrices and so that $\det M_G=\pm (n+1)$ and $\det M_H=\pm (m+1)$.

\medskip
Denote the matrix in the array in the proof of Lemma~\ref{Lem 4.1}
by $M_K$. It is clear that $\det M_K=\pm |K|$. The extended matrices $\ti{M}_G$, $\ti{M}_H$, $\ti{M}_K$ are
\[\ti{M}_G=\left(\begin{array}{c}
M_G\\
s\end{array}\right),\quad
\ti{M}_H=\left(\begin{array}{c}
M_H\\
t\end{array}\right), \quad
\ti{M}_K=\left(\begin{array}{ccc}
M_G & 0 & X\\
0 & M_H & Y\\
0 & 0 & I\\
s & t & e \end{array}\right),\]
where $s=(s_1,\ldots, s_n)$ and $t=(t_1,\ldots,t_m)$ have entries from the set $\{-1,0,1\}$, $e=(1,\ldots,1)$, $I$ is the $nm \times nm$
identity matrix, and $X,Y$ are the two blocks we also see in the array in the proof of Lemma~\ref{Lem 4.1}.
We have to show that $\ti{M}_K$ is a basis matrix for $\cL(K)$, and since $|K|=(n+1)(m+1)$, this
is equivalent to the equality $\det \ti{M}_K^\top \ti{M}_K=(n+1)^3(m+1)^3$.

\medskip
We expand $\det \ti{M}_K^\top \ti{M}_K$ by the Cauchy-Binet formula.
In what follows we also write $|A|$ for the determinant of a matrix $A$.
We then have
\[|\ti{M}_K^\top \ti{M}_K|=|M_\ell|^2+\sum_{j,k}|M_{j,k}|^2+\sum_k |{M}_{H,k}|^2+\sum_j |{M}_{G,j}|^2,\]
the matrices on the right resulting from $\ti{M}_K$ after deleting the last row, the row labelled by $(g_j,h_k)$,
the row labelled by $(0,h_k)$, and the row labelled by $(g_j,0)$, respectively. The matrix $M_\ell$ is
upper block-triangular and hence $|M_\ell|^2=|M_G|^2|M_H|^2$. We may expand the determinant $|M_{jk}|$
along the rows intersecting the identity matrix $I$, giving
\[|M_{jk}|^2=\left|\begin{array}{ccc}
M_G & 0 & X_j\\
0 & M_H & Y_k\\
s & t & 1\end{array}\right|^2,\]
where $X_j$ and $Y_k$ are columns with a single $-1$ and zeros otherwise. Adding the first $n+m$ rows to the last row,
we get
\[|M_{jk}|^2=\left|\begin{array}{ccc}
M_G & 0 & X_j\\
0 & M_H & Y_k\\
0 & 0 & -1\end{array}\right|^2=|M_G|^2|M_H|^2.\]
Expanding the determinant $|{M}_{H,k}|$ along the rows which intersect the identity matrix $I$ we obtain
\[|{M}_{H,k}|^2=\left|\begin{array}{cc}
M_G & 0\\
p & \ti{M}_{H,k}\end{array}\right|^2, \quad p=\left(\begin{array}{c}0\\ s\end{array}\right), \quad
\ti{M}_{H,k}=\left(\begin{array}{c} M_{H,k}\\ t\end{array}\right),\]
where $M_{H,k}$ arises from $M_H$ by deleting the $k$th row. The matrix $\ti{M}_{H,k}$ is square and hence
\[|{M}_{H,k}|^2=|M_G|^2 |\ti{M}_{H,k}|^2.\]
Analogously, $|{M}_{G,j}|^2=|\ti{M}_{G,j}|^2|M_H|^2$. In summary,
\[|\ti{M}_K^\top \ti{M}_K|=|M_G|^2|M_H|^2+nm |M_G|^2|M_H|^2+|M_G|^2 \sum_k |\ti{M}_{H,k}|^2+\sum_j
|M_H|^2 |\ti{M}_{G,j}|^2.\]
Again due to Cauchy-Binet,
\begin{eqnarray*}
& & |M_H|^2+\sum_k |\ti{M}_{H,k}|^2=|\ti{M}_H^\top \ti{M}_H|=(m+1)^3, \\
& & |M_G|^2+\sum_j |\ti{M}_{G,j}|^2=|\ti{M}_G^\top \ti{M}_G|=(n+1)^3,
\end{eqnarray*}
and taking into account that $|M_G|^2=(n+1)^2$ and $|M_H|^2=(m+1)^2$, we arrive at
the conclusion that $|\ti{M}_K^\top \ti{M}_K|$ is equal to
\[(n+1)^2(m+1)^2(1+nm)+(n+1)^2[(m+1)^3-(m+1)^2]
+(m+1)^2[(n+1)^3-(n+1)^2],\]
which equals $(n+1)^3(m+1)^3$, as desired. $\;\:\square$

\medskip
In Section \ref{SC} we showed that the hypothesis of Lemma~\ref{Lem 4.1a} is satisfied if
$G$ and $H$ are cyclic groups of order at least five. Successive application of Lemma~\ref{Lem 4.1a}
therefore gives Theorem~\ref{Theo 1.3} for all groups $\bZ_{m_1} \times \cdots \times \bZ_{m_k}$
with $m_1, \ldots, m_k \ge 5$.

\begin{lem} \label{Lem 4.2a}
Let $m \in \{2,3,4\}$ and let $G$ be a finite Abelian group such that $\cL(G)$ has a basis of minimal vectors
and such that $d(G)=2$. Also suppose that there is an admissible array $\bg||M_G$
coming from minimal basis vectors such that $\det M_G=\pm |G|$. Put $K=\bZ_m \times G$.
Then $\cL(K)$ has a basis of minimal vectors, $d(K)=2$,
and there exists an admissible array $\bk||M_K$ resulting from minimal basis vectors such that $\det M_K=\pm |K|$.
\end{lem}

{\em Proof.} We proceed as in the proof of the preceding lemma.
Let $G$, $\bg$, and the admissible arrays $\bk||M_K$ be as in the proof of Lemma~\ref{Lem 4.2}.
These arrays are associated with vectors of minimum length $2$ and the extended matrices $\ti{M}_K$ are of the form
\[\ti{M}_K =\left( \begin{array}{ccc}
M_G & * & * \\
0 & M_m & * \\
0 & 0 & I \\
s & t & e \end{array} \right).\]
We already know that $\det M_G=\pm |G| = \pm(n+1)$ and $\det M_m=\pm m$. It remains to prove that
$\det \ti{M}_K^\top \ti{M}_K=m^3(n+1)^3$.

\medskip
We consider the case $m=3$. The cases $m=2$ and $m=4$ may be disposed of in a completely analogous fashion.
Expanding $\det \ti{M}_K^\top \ti{M}_K$ by the Cauchy-Binet formula we get
\[\det \ti{M}_K^\top \ti{M}_K=|M_\ell|^2 +\sum_{k=1}^{2n-1}|M_{I,k}|^2+\sum_{j=1}^3 |M_{3,j}|^2+\sum_{i=1}^n|M_{G,i}|^2,\]
where $M_\ell$ results from deleting the last row, $M_{I,k}$ comes from deleting the row which contains the $k$th
entry $1$ of
the $(2n-1) \times (2n-1)$ identity matrix $I$, $M_{3,j}$ arises from deleting the row containing the $j$th row of $M_3$,
and $M_{G,i}$ emerges  from deleting the $i$th row. Clearly, $|M_\ell|^2=|M_G|^2|M_3|^2=9(n+1)^2$. Expanding $|M_{I,k}|$
along the rows intersecting the identity matrix and  adding after that the first $n+3$ rows to the last row, we obtain
\[|M_{I,k}|^2=\left|\begin{array}{ccc}
M_G & * & *\\
0 & M_3 & *\\
s & t & 1\end{array}\right|^2=\left|\begin{array}{ccc}
M_G & * & *\\
0 & M_3 & *\\
0 & 0 & -1\end{array}\right|^2=|M_G|^2|M_3|^2=9(n+1)^2.\]
We expand $M_{3,j}$ again along the rows intersecting the identity matrix and then add the first $n+2$ rows to the last.
What results is
\[|M_{3,j}|^2=\left|\begin{array}{cc}
M_G & * \\
0 & Q_{3,j}\\
s & t\end{array}\right|^2=
\left|\begin{array}{cc}
M_G & * \\
0 & Q_{3,j}\\
0 & t_j\end{array}\right|^2=|M_G|^2 \left|\begin{array}{c} Q_{3,j}\\t_j\end{array}\right|^2=(n+1)^2\left|\begin{array}{c} Q_{3,j}\\t_j\end{array}\right|^2\]
with
\[\sum_{j=1}^3\left|\begin{array}{c} Q_{3,j}\\t_j\end{array}\right|^2=
\left|\begin{array}{rrr}1 & 1 & 1\\-1 & 1 & 0\\-1 & 0 & 1\end{array}\right|^2
+\left|\begin{array}{rrr}1 & 0 & -1\\-1 & 1 & 0\\-1 & -1 & -1\end{array}\right|^2
+\left|\begin{array}{rrr}1 & 0 & -1\\1 & 1 & 1\\1 & -1 & 0\end{array}\right|^2.\]
Each determinant on the right equals $3$ and hence the sum of their squares is $27$. Finally,
again after expansion along the rows intersecting $I$ and a row change,
\[|M_{G,i}|^2=\left|\begin{array}{cc}\ti{M}_{G,i} & * \\0 & M_3\end{array}\right|^2
=|\ti{M}_{G,i}|^2|M_3|^2=9|\ti{M}_{G,i}|^2.\]
By the Cauchy-Binet formula,
\[\sum_{j=1}^n |\ti{M}_{G,i}|^2=|\ti{M}_G^\top \ti{M}_G|-|M_G|^2=(n+1)^3-(n+1)^2=(n+1)^2n.\]
Putting things together we see that
\[\det\ti{M}_K^\top\ti{M}_K=9(n+1)^2(1+2n-1)+27(n+1)^2+9(n+1)^2n=3^3(n+1)^3,\]
which is what we wanted. $\;\:\square$

\begin{lem} \label{Lem 4.3a}
Let $G$ be one of the groups $\bZ_2 \times \bZ_4$, $\bZ_3 \times \bZ_3$, $\bZ_4 \times \bZ_4$. Then $\cL(G)$
has a basis of minimal vectors, $d(G)=2$, and there exists an admissible array $\bg||M_G$
coming from minimal basis vectors such that $\det M=\pm |G|$ and $\det \ti{M}_G^\top \ti{M}_G=|G|^3$.
\end{lem}

{\em Proof.}
The admissible array
$\bg||M_7$ we have shown for $G=\bZ_2\times \bZ_4$ in the proof of Lemma~\ref{Lem 4.3}
satisfies $\det M_7=8$ and $\det \ti{M}_7^\top \ti{M}_7 =8^3$.
The array
\[\bg||M_8=\begin{array}{l||rrrr|rrrr}
(0,1)& 1 &   &   &  &-1 &   & -1&  \\
(0,2)&   & 1 & 1 &  &   & -1&   & \\
(1,0)& 1 &   &   & 1& 1 & 1 &   & \\
(1,1)&-1 &   &   & 1&   &   &   & \\
\hline
(1,2)&   &   &   &  &   &   & 1 & -1\\
(2,0)&   & 1 & -1&  &   &   &   & 1 \\
(2,1)&   &   & 1 &-1& 1 &   &   & \\
(2,2)&   &-1 &   &  &   & 1 & 1 & 1\\
\end{array}\]
is admissible for $\bZ_3 \times \bZ_3$, and we have $\det M_8=9$ and $\det \ti{M}_8^\top\ti{M}_8=9^3$.
The array $\bg||M_{15}$ given by
\[\begin{array}{l|rrr|rrrr|rrrr|rrrr}
(0,1) & 1 & 1 & 1 & 1 &   &   &   &   &   &   &   &   &   &   &  \\
(0,2) &   &   &   &   & 1 &   &   &   &   &   &   &   &   &   &  \\
(0,3) &   &   &   &   &   & 1 & 1 & 1 &   &   &   &   &   &   &  \\
\hline
(1,0) &   & -1&   &   & 1 &-1 &   &   &   &   &   &   &   &   &  \\
(1,1) &   &   &   &   &   & 1 &   &   &   &   &   &   &   &   &  \\
(1,2) & 1 &   &   &   &-1 &   &   &   &   &   &   &-1 &   &-1 &  \\
(1,3) &-1 & 1 &   &   &   &   &   &   & 1 & 1 & 1 &   &-1 &   &  \\
\hline
(2,0) &   &   &   &-1 &   &   &   &   &   &   &   & 1 & 1 &   &  \\
(2,1) &   &   &   &   &   &   &   &   & 1 &   &   &   &   & 1 &  \\
(2,2) &   &   & 1 &   &   &   &   &   &   & 1 &   &   &   &   &  \\
(2,3) &   &   & -1& 1 &   &   &   &   &   &   & 1 &   &   &   &-1 \\
\hline
(3,0) &   &   &   &   &   &   &-1 &   &-1 &   &   &   &   &   &  \\
(3,1) &   &   &   &   &   &   & 1 &-1 &   &-1 &   &   &   & 1 & 1 \\
(3,2) &   &   &   &   &   &   &   & 1 &   &   &-1 & 1 &   &   & 1 \\
(3,3) &   &   &   &   &   &   &   &   &   &   &   &   & 1 &   & \end{array}\]
is admissible for $\bZ_4 \times \bZ_4$ with $\det M_{15}=-16$ and $\det \ti{M}_{15}^\top\ti{M}_{15}=16^3$.
$\:\;\square$

\medskip
The proof of Theorem \ref{Theo 1.3} may now be completed as at the end of Section~\ref{SD}.

\section{Bounds for the covering radius} \label{SCov}

In this section we study the lattices $\cL_n:=\cL(\bZ_m)$ ($m=n+1$) with the goal of proving
Theorem~\ref{Theo 1.4}. For $n \ge 2$, let $B_{n+1,n}$ be the $(n+1)\times n$ version of the matrices
\[B_{3,2}=\left(\begin{array}{rr}
-2 & 1\\
1 & -2\\
1 & 1\end{array}\right), \;
B_{4,3}=\left(\begin{array}{rrr}
-2 & 1 & 0\\
1 & -2 & 1\\
0 & 1 & -2\\
1 & 0 & 1\end{array}\right), \;
B_{5,4}=\left(\begin{array}{rrrr}
-2 & 1 & 0 & 0 \\
1 & -2 & 1 & 0 \\
0 & 1 & -2 & 1 \\
0 & 0 & 1 & -2 \\
1 & 0 & 0 & 1 \\
\end{array}\right).\]
Note that in Section~\ref{SC} we denoted these matrices by $\ti{T}_{n+1}$. In other words, we now denote $\ti{T}_m$
by $B_{m,m-1}$. The $(n+1) \times k$ matrix formed by
first $k$ columns of $B_{n+1,n}$ is denoted by $B_{n+1,k}$.

\begin{exa} \label{Exa 9.1}
{\rm Let us begin with an example. Consider $G=\bZ_4$. We know from Section~\ref{SC}
that
\[B_{4,3} \; (=\ti{T}_4)\; =\left(\begin{array}{rrr}
-2 & 1 & 0\\
1 & -2 & 1\\
0 & 1 & -2\\
\hline
1 & 0 & 1
\end{array}\right)\]
is a basis matrix for the lattice $\cL_3:=\cL(G)$.
This follows from the fact that
\[V_3:=\sqrt{\det \ti{B}_3^\top \ti{B}_3} = 8= 4^{3/2}.\]

Let ${\bob}_1,{\bob}_2,{\bob}_3$ be the columns of ${B}_{4,3}$.
Then
\[{B}_{4,2}=(\bob_1 \;\bob_2)=\left(\begin{array}{rr}
-2 & 1 \\
1 & -2 \\
0 & 1 \\
\hline
1 & 0
\end{array}\right), \quad
{B}_{4,1}=(\bob_1)=\left(\begin{array}{r}
-2 \\
1 \\
0 \\
\hline
1
\end{array}\right).
\]
Let further $\cL_2$ and $\cL_1$ be the sublattices of $\cL_3$ spanned by the columns of ${B}_{4,2}$ and ${B}_{4,1}$. The determinants of $\cL_2$
and $\cL_1$ are
\[V_2=\sqrt{\det {B}_{4,2}^\top {B}_{4,2}}=\left|\begin{array}{rr}
6 & -4\\
-4 & 6 \end{array}\right|^{1/2}=\sqrt{20}, \quad V_1=\sqrt{\det{B}_{4,1}^\top {B}_{4,1}}=\sqrt{6}.\]
The lattice $\cL_1$
is spanned by a vector of length $\sqrt{6}$ and can therefore be covered by $1$-dimensional balls of radius $r_1=\sqrt{6}/2$ centered at the
lattice points.
Now consider an arbitrary point $x$ in ${\rm span}_\bR \cL_2$. We may assume that this point lies between
the two lines ${\rm span}_\bR \cL_1$ and ${\bob}_2 +{\rm span}_\bR \cL_1$. Let $h_1$ be the distance between these two lines.
The distance between $x$ and one of the two lines is at most $h_1/2$. This implies that $x$ is contained in a $2$-dimensional ball
of radius $r_2\le\sqrt{r_1^2+(h_1/2)^2}$ centered at a lattice point of $\cL_1$ or ${\bob}_2+\cL_1$.
Since the area of a parallelogram is the product of the length of the baseline and the height,
we have $V_2=V_1h_1$. Thus, ${\rm span}_\bR\cL_2$ may be covered by $2$-dimensional balls of radius $r_2$ centered at the points
of $\cL_2$ where,
\[r_2^2\le r_1^2+\left(\frac{V_2}{2V_1}\right)^2=\frac{6}{4}+\frac{20}{4\cdot 6}=\frac{23}{12}.\]
Now take a point $y$ in ${\rm span}_\bR\cL_3$, without loss of generality between the two planes
${\rm span}_\bR\cL_2$ and ${\bob}_3+{\rm span}_\bR \cL_2$. Letting $h_2$ denote the distance between these two planes,
there is a point in $\cL_2$ or $\bob_3+\cL_2$ whose distance to $y$ is at most $\sqrt{r_2^2+(h_2/2)^2}$. Since $V_3=V_2h_2$,
we conclude that ${\rm span}_\bR\cL_3$ may be covered by $3$-dimensional balls of radius $r_3$ with the centers at
the points of $\cL_3$, where
\[r_3^2 \le r_2^2+\left(\frac{V_3}{2V_2}\right)^2 \le \frac{23}{12}+\frac{64}{4\cdot 20}=\frac{47}{15}.\]
Consequently, $\mu(\bZ_4)\le \sqrt{47/15} \approx 1.7701$.
}\end{exa}

\begin{prop} \label{Prop 9.2}
Let $\bob_1, \ldots, \bob_n$ be points in the root lattice $\cA_n$ such that
\[{\rm span}_\bR\{\bob_1, \ldots, \bob_n\}={\rm span}_\bR\cA_n.\]
For $k=1, \ldots, n$, denote by $\cL_k$ the lattice spanned by $\bob_1, \ldots, \bob_k$,
let $C_{n+1,k}$ stand for the $(n+1) \times k$ matrix whose columns are the coordinates of $\bob_1, \ldots, \bob_k$
and put \[V_k =\sqrt{\det C_{n+1,k}^\top C_{n+1,k}}.\] If ${\rm span}_\bR\cL_k$ {\rm (}$1 \le k \le n-1${\rm )} can be covered
by $k$-dimensional balls of radius $r_k$ centered at the points of $\cL_k$, then ${\rm span}_\bR\cL_{k+1}$ can be covered
by balls of dimension $k+1$ centered at the points of $\cL_{k+1}$ whose radius $r_{k+1}$ satisfies
\[r_{k+1}^2 \le r_k^2+\left(\frac{V_{k+1}}{2V_k}\right)^2,\]
and consequently,
\[r_n^2 \le r_1^2+\left(\frac{V_{2}}{2V_1}\right)^2+ \cdots + \left(\frac{V_{n}}{2V_{n-1}}\right)^2.\]
\end{prop}

{\em Proof.} This can be shown by the argument employed in Example~\ref{Exa 9.1}. $\;\:\square$

\bigskip
The only problem in general is the computation of the determinants $V_k$. Fortunately, this is easy for $\cL_n=\cL(\bZ_m)$
($m=n+1$), in which case the matrices $C_{n+1,k}$ are just the matrices $B_{n+1,k}$ we introduced above.
We also need the $n \times n$ versions $Q_n$ of the matrices
\[Q_2=\left(\begin{array}{rr}
6 & -3\\
-3 & 6\end{array}\right),\;
Q_3=\left(\begin{array}{rrr}
6 & -4 & 2\\
-4 & 6 & -4\\
2 & -4 & 6\end{array}\right),\;
Q_4=\left(\begin{array}{rrrr}
6 & -4 & 1 & 1\\
-4 & 6 & -4 & 1\\
1 & -4 & 6 & -4\\
1 & 1 & -4 & 6 \end{array}\right),\]
\[
Q_5= \left(\begin{array}{rrrrr}
6 & -4 & 1 & 0  & 1\\
-4 & 6 & -4 & 1  & 0\\
1 & -4 & 6 & -4 & 1\\
0 & 1 & -4 & 6 & -4\\
1  & 0 & 1 & -4 & 6
\end{array}\right),\;
Q_6= \left(\begin{array}{rrrrrr}
6 & -4 & 1 & 0 & 0 & 1\\
-4 & 6 & -4 & 1 & 0 & 0\\
1 & -4 & 6 & -4 & 1 & 0\\
0 & 1 & -4 & 6 & -4 & 1\\
0 & 0 & 1 & -4 & 6 & -4\\
1 & 0 & 0 & 1 & -4 & 6
\end{array}\right).\]
Finally, for $k \ge 1$, we denote by $R_k$ the $k \times k$ version of the matrices $R_1=(6)$,
\[
R_2=\left(\begin{array}{rr}
6 & -4\\
-4 & 6\end{array}\right),\;
R_3=\left(\begin{array}{rrr}
6 & -4 & 1\\
-4 & 6 & -4\\
1 & -4 & 6\end{array}\right),\;
R_4=\left(\begin{array}{rrrr}
6 & -4 & 1 & 0\\
-4 & 6 & -4 & 1\\
1 & -4 & 6 & -4\\
0 & 1 & -4 & 6 \end{array}\right),\]
\[
R_5= \left(\begin{array}{rrrrr}
6 & -4 & 1 & 0  & 0\\
-4 & 6 & -4 & 1  & 0\\
1 & -4 & 6 & -4 & 1\\
0 & 1 & -4 & 6 & -4\\
0 & 0 & 1 & -4 & 6
\end{array}\right),\;
R_6= \left(\begin{array}{rrrrrr}
6 & -4 & 1 & 0 & 0 & 0\\
-4 & 6 & -4 & 1 & 0 & 0\\
1 & -4 & 6 & -4 & 1 & 0\\
0 & 1 & -4 & 6 & -4 & 1\\
0 & 0 & 1 & -4 & 6 & -4\\
0 & 0 & 0 & 1 & -4 & 6
\end{array}\right).\]

\begin{lem} \label{Lem 9.3}
For $n \ge 2$ and $1 \le k \le n-1$,
\begin{eqnarray*}
& & B_{n+1,n}^\top B_{n+1,n}=Q_n, \quad \det Q_n=(n+1)^3,\\
& & B_{n+1,k}^\top B_{n+1,k}=R_k, \quad \det R_k=\frac{(k+1)(k+2)^2(k+3)}{12}.
\end{eqnarray*}
\end{lem}

{\em Proof.} The formulas for the products of the matrices can be verified by straightforward
computation. The formula for $\det Q_n$ is nothing but Lemma~\ref{Lem Cyc2}.
The formula for $\det R_k$ was first established in~\cite{BoSiJFA}. Proofs of that formula can also be found in~\cite[Theorem~10.59]{BoSiBu} and~\cite{BoWi}.
$\;\:\square$

\bigskip
{\em Proof of Theorem \ref{Theo 1.4}.} We know from Section \ref{SC} that $B_{n+1,n}=(\bob_1, \ldots,\bob_n)$ is a basis matrix for $\cL_n:=\cL(\bZ_{n+1})$.
Let $\cL_k$ be the sublattices as in Proposition~\ref{Prop 9.2}. The $1$-dimensional lattice $\cL_1$ is spanned by a vector of length $\sqrt{6}$.
We may therefore use Proposition~\ref{Prop 9.2} with $r_1=\sqrt{6}/2$ to obtain that
\[\mu(\cL_n)^2  \le  \frac{6}{4}+\left(\frac{V_{2}}{2V_1}\right)^2+ \cdots + \left(\frac{V_{n}}{2V_{n-1}}\right)^2
 =  \frac{6}{4}+\frac{1}{4}\sum_{k=1}^{n-2}\frac{V_{k+1}^2}{V_k^2}+\frac{1}{4}\,\frac{V_n^2}{V_{n-1}^2}.\]
From Lemma \ref{Lem 9.3} we see that if $1 \le k \le n-2$, then
\begin{eqnarray*}
\frac{V_{k+1}^2}{V_k^2} & = & \frac{\det B_{n+1,{k+1}}^\top B_{n+1,{k+1}}}{\det B_{n+1,k}^\top B_{n+1,k}}
= \frac{(k+2)(k+3)^2(k+4)}{(k+1)(k+2)^2(k+3)}\\
& = & \frac{(k+3)(k+4)}{(k+1)(k+2)}=1+\frac{2(2k+5)}{(k+1)(k+2)}=1+2\left(\frac{3}{k+1}-\frac{1}{k+2}\right).
\end{eqnarray*}
Consequently,
\begin{eqnarray*}
\sum_{k=1}^{n-2}\frac{V_{k+1}^2}{V_k^2} & = & n-2+2\sum_{k=1}^{n-2}\left(\frac{3}{k+1}-\frac{1}{k+2}\right)
 =  n-2+2\left(\frac{3}{2}+2\sum_{k=3}^{n-1}\frac{1}{k}-\frac{1}{n}\right)\\
 & = & n+1-\frac{2}{n}+4\sum_{k=3}^{n-1}\frac{1}{k} < n+1-\frac{2}{n}+4\int_2^{n-1}\frac{dx}{x}\\
& = & n+1-\frac{2}{n}+4\log(n-1)-4\log 2.
\end{eqnarray*}
Lemma \ref{Lem 9.3} also implies that
\[\frac{V_n^2}{V_{n-1}^2} =\frac{\det B_{n+1,n}^\top B_{n+1,n}}{\det B_{n+1,n-1}^\top B_{n+1,n-1}}
=\frac{12(n+1)^3}{n(n+1)^2(n+2)}=\frac{12(n+1)}{n(n+2)}.\]
In summary,
\begin{eqnarray*}
\mu(\cL_n)^2 & < & \frac{1}{4}\left(6+n+1-\frac{2}{n}+4\log(n-1)-4\log 2+\frac{12(n+1)}{n(n+2)}\right)\\
& = &  \frac{1}{4}\left(n +4\log (n-1) + 7-4\log 2+ \frac{10n+8}{n(n+2)}\right),
\end{eqnarray*}
and since $(10n+8)/(n(n+2)) < 10/n$, we arrive at the asserted bound.
$\;\:\square$

\bigskip
We remark that Example \ref{Exa 9.1} gives $\mu(\bZ_4) < 1.7701$ whereas the table presented in Section~\ref{SI} shows the slightly
worse bound $\mu(\bZ_4) < 1.8257$. This discrepancy is caused by the circumstance that in
Example~\ref{Exa 9.1} we didn't estimate a sum by an integral.

\medskip
As already mentioned in the introduction, Sha~\cite{Sha} showed that $\mu(G) \le \mu(\cA_n)+\sqrt{2}$
if $G$ is a group coming from elliptic curves over finite fields. Actually, his proof works for
arbitrary finite Abelian groups $G$. It goes as follows. Let
$\xi=(\xi_1, \ldots, \xi_n,\xi_0)  \in {\rm span}_\bR \cA_n$,
where $\xi_0:=-\xi_1-\cdots-\xi_n$, and pick $v=(v_1,\ldots,v_n, v_0) \in \cA_n$  as a point
for which $\n \xi-v\n \le \mu(\cA_n)$. Then one may proceed as in~\cite[proof of Theorem 3.4]{FM}. Namely, let
$v_1g_1+\cdots+v_ng_n=g_j$. If $g_j$ is not the zero of the group, put
\[x=(x_1, \ldots, x_n, x_0)=(v_1, \ldots, v_{j-1}, v_j-1, v_{j+1}, \ldots, v_n, v_0+1).\]
Then $x \in \cL(G)$ and $\n v-x\n =\sqrt{2}$. In case $g_j$ is the zero of the group, let $x=v$, so that
$x \in \cL(G)$ and $\n v-x\n=0$. In either case, $\n\xi-x\n \le \mu(\cA_n)+\sqrt{2}$.

\medskip
The only difference between \cite{FM} and \cite{Sha} is that in \cite{FM}
the point $v=(v_1,\ldots,v_n, v_0) \in \cA_n$ was chosen so that
$v_i$ is the nearest integer to $\xi_i$ for $i=1,\ldots, n$. To ensure that $v$
is in $\cA_n$, one had to take $v_0=-v_1-\cdots-v_n$, and as the difference between $\xi_0$ and $v_0$
may be large, the bound for the covering radius obtained in~\cite{FM} was too coarse.
Sha's clever choice of $v=(v_1,\ldots,v_n, v_0) \in \cA_n$  as a point
for which $\n \xi-v\n \le \mu(\cA_n)$ remedied this defect.

\section{The automorphism group}\label{aut}

In this section we start out with the proof of Theorem~\ref{Theo 1.5}.
\medskip

{\em Proof of Theorem~\ref{Theo 1.5}.} Let $G = \left\{ 0,g_1,\dots,g_n \right\}$ be a finite Abelian group and recall that
$$\cL(G) = \left\{ X = \left(x_1,\dots,x_n, -\sum_{i=1}^n x_i \right) \in \bZ^{n+1} : \sum_{i=1}^n x_i g_i = 0 \right\}.$$
Every automorphism of $G$ fixes $0$ and permutes the elements $g_1,\dots,g_n$. Hence $\Aut(G)$ can be identified (via a canonical isomorphism) with a
subgroup of the symmetric group~$S_n$. We denote this subgroup by $H$. Our objective is to construct a group isomorphism $\Phi: H \to
\Aut(\cL(G)) \cap S_n$, where $\Aut(\cL(G))$ on the right is identified with a subgroup of $\GL_n(\bZ)$
as described in Section~\ref{SI} and $S_n$ on the right is viewed in the natural fashion as the subgroup of the permutation matrices in~$\GL_n(\bZ)$.
\medskip

Let $\si \in H$. Then, for every $g_i \in G$, $\si(g_i) = g_{\si(i)}$ and $\si(0)=0$. If
$$X = \left(x_1,\dots,x_n, -\sum_{i=1}^n x_i \right) \in \cL(G),$$
then $\sum_{i=1}^n x_i g_i = 0$. Notice that $\si^{-1}$ is also in $H$, and so
$$0 = \si^{-1}(0) = \sum_{i=1}^n x_i g_{\si^{-1}(i)} = \sum_{i=1}^n x_{\si(i)} g_i.$$
Now define $\tau=\Phi(\si)$ on $\cL(G)$ by
$$\tau \left(x_1,\dots,x_n, -\sum_{i=1}^n x_i \right) := \left(x_{\si(1)},\dots,x_{\si(n)}, -\sum_{i=1}^n x_{\si(i)} \right).$$
It is clear that $\tau$ maps $\cL(G)$ onto itself.
The matrix $U \in \GL_n(\bZ)$ corresponding to $\tau$ as described in Section~\ref{SI} is obviously a permutation matrix.
Consequently, $\tau$ is in $\Aut(\cL(G)) \cap S_n$. Finally, it is readily seen that $\Phi$ is an injective group
homomorphism. Hence $\Phi(H) \leq \Aut(\cL(G)) \cap S_n$.
\medskip

It remains to show that $\Phi(H) = \Aut(\cL(G)) \cap S_n$. So suppose $\tau \in \Aut(\cL(G)) \cap S_n$. If
$$X = \left( x_1,\dots,x_n,-\sum_{i=1}^n x_i \right) \in \cL(G),$$
then $\tau(X) = (x_{\si(1)},\dots,x_{\si(n)},-\sum_{i=1}^n x_{\si(i)})$ with some $\si \in S_n$, and since
both $X$ and $\tau(X)$ belong to  $\cL(G)$, it follows that
$$0 = \sum_{i=1}^n x_i g_i = \sum_{i=1}^n x_{\si(i)} g_i.$$
We have $\tau =\Phi(\si)$ with $\si: G \to G$ defined by
$\si(g_i) := g_{\si(i)}$ and $\si(0) := 0$. To complete the proof, we only need to show that $\si$ is a group homomorphism,
i.e. that
$$\si(g_i + g_j) = g_{\si(i)} + g_{\si(j)}.$$
Since $g_i + g_j \in G$, there must be some $g_k \in G$ such that $g_i + g_j = g_k$. In other words
$$g_i + g_j - g_k = 0.$$
Therefore the vector $X$ with $i$th and $j$th coordinates equal to $1$, $k$th coordinate equal to $-1$, $(n+1)$st coordinate equal to $- (1 + 1 - 1) = -1$,
and the rest of the coordinates equal to $0$ must be in $\cL(G)$. Hence the vector $\tau(X)$ also lies in $\cL(G)$. This vector has
 $\si(i)$th and $\si(j)$th coordinates equal to $1$,
$\si(k)$th coordinate equal to $-1$, $(n+1)$st coordinate equal to $-1$, and the rest of the coordinate equal to 0. This means that the equality
$$g_{\si(i)} + g_{\si(j)} - g_{\si(k)} = 0$$
must be satisfied in $G$, and hence
$$\si(g_i + g_j) = \si(g_k) = g_{\si(k)} = g_{\si(i)} + g_{\si(j)}.$$
In summary, $\si \in H$, and so $\Aut(\cL(G)) \cap S_n =\Phi(H)$, as desired.
$\;\:\square$

\bigskip
Theorem \ref{Theo 1.5} has an interesting geometric interpretation in terms of the theory of quadratic forms (see, for instance,~\cite{schur_book} for a detailed account of this subject and its connections to lattice theory). A real quadratic form in $n$ variables $X = (x_1,\dots,x_n)^\top$ can always be written in a unique way as
$$q(X) = X^\top A X,$$
where $A$ is an $n \times n$ real symmetric matrix. Hence the space of real quadratic forms in $n$ variables can be identified with the space of their coefficient matrices, which is the $\binom{n+1}{2}$-dimensional real vector space $\s^n$ of $n \times n$ real symmetric matrices. The set of positive definite forms $\s^n_{>0}$ is an open convex cone in $\s^n$ given by $n$ polynomial inequalities (the Sylvester criterion).

\medskip
Let $B$ be an $m \times n$ real matrix of rank $n$, $1 \leq n \leq m$, then $\cL = B\bZ^n$ is a lattice of rank $n$ in $\R^m$.
The so-called norm form of $\cL$, corresponding to the choice of the basis matrix $B$,
is defined as the positive definite quadratic form in $n$ variables, given by
$$q_B(X) = X^\top (B^\top B) X.$$
The function $B \mapsto B^\top B$ induces a bijection between the space of lattices (up to isometry) and the cone $\s^n_{>0}$ of positive definite quadratic forms (up to arithmetic equivalence).

\medskip
Given a form $q \in \s^n$, its automorphism group is defined by
$$\Aut(q) := \left\{ \tau \in \GL_n(\bZ) : q(\tau(X)) = q(X)\;\: \mbox{for all} \;\: X \in \R^n \right\}.$$
This is a finite group: it is contained in the intersection of the discrete group $\GL_n(\bZ)$ with the compact group $\cO_n(\R)$ of real orthogonal matrices. If $q \in \s^n_{>0}$ and $\cL$ is the corresponding lattice, then $\Aut(q) = \Aut(\cL)$, the automorphism group of $\cL$. Furthermore, given any finite subgroup $H$ of $\GL_n(\bZ)$, there exists
a $q \in \s_{>0}^n$ (and hence a lattice) with $H \leq \Aut(q)$. Indeed, if $H \leq \GL_n(\bZ)$ is a finite group and $f \in \s^n_{>0}$, then
the form defined by
$$q(X) := \sum_{\tau \in H} f(\tau(X))$$
is in $\s^n_{>0}$ and $H \leq \Aut(q)$. Finally, for a fixed finite subgroup $H$ of $\GL_n(\bZ)$,  define
$$\cB(H) = \left\{ q \in \s^n : H \leq \Aut(q) \right\}. $$
The set $\cB(H)$ is not empty by the above remark, and hence it is easily seen to be a real vector space. It is called the Bravais manifold of $H$. Define also the open convex polyhedral cone $\cB_{>0}(H) = \cB(H) \cap \s^n_{>0}$, which can be identified with the set of all lattices whose automorphism groups contain $H$.
\medskip

Investigation of properties of Bravais manifolds corresponding to different finite subgroups of~$\GL_n(\bZ)$ is of interest in lattice theory. Our Theorem~\ref{Theo 1.5} implies that the lattice $\cL(G)$ coming from an Abelian group $G$ of order $n+1$ via our construction is contained in the Bravais cone $\cB_{>0}(\Aut(G))$.

\bigskip
{\bf Acknowledgement.} We thank Min Sha for kindly informing us of his results prior to posting the preprint \cite{Sha}.

\bigskip
A. B\"ottcher, Fakult\"at f\"ur Mathematik, TU Chemnitz, 09107 Chemnitz, Germany

{\tt aboettch@mathematik.tu-chemnitz.de}

\medskip
L. Fukshansky, Department of Mathematics,  Claremont McKenna College,

850 Columbia Ave,
Claremont, CA 91711, USA

{\tt lenny@cmc.edu}

\medskip
S. R. Garcia, Department of Mathematics, Pomona College,

610 N. College Ave, Claremont, CA 91711, USA

{\tt stephan.garcia@pomona.edu}, URL: \url{http://pages.pomona.edu/~sg064747/}

%URL: {\tt http://pages.pomona.edu/$\scriptstyle\mathtt{\sim}$sg064747/}

\medskip
H. Maharaj, 8543 Hillside Road, Rancho Cucamonga, CA 91701, USA

{\tt hmahara@g.clemson.edu}
\end{document}